\documentclass[a4paper, 11pt, twoside]{article}
\usepackage{amsmath}
\usepackage{amscd}
\usepackage{amssymb}
\usepackage{amsthm}
\usepackage{graphics}
\usepackage{indentfirst}
\usepackage{latexsym}
\usepackage{amsfonts}
\usepackage{multicol}
\usepackage{stmaryrd}
\usepackage{array}
\usepackage{youngtab}
\usepackage{pxfonts}
\usepackage{verbatim}
\usepackage{bbm}
\usepackage{pgf}
\usepackage{tikz}
\usepackage{nccmath}

\makeatletter
\newenvironment{pf}[1][\proofname] {\par\pushQED{\qed}\normalfont\topsep6\p@\@plus6\p@\relax\trivlist\item[\hskip\labelsep\bfseries#1\@addpunct{.}]\ignorespaces}{\popQED\endtrivlist\@endpefalse}
\makeatother

\newtheoremstyle{mattthm}{}{}{\itshape}{}{\bfseries}{.}{ }{}

\theoremstyle{mattthm}
\newtheorem{lemma}{Lemma}[section]
\newtheorem{propn}[lemma]{Proposition}
\newtheorem{thm}[lemma]{Theorem}
\newtheorem{cory}[lemma]{Corollary}

\newtheoremstyle{mattdef}{}{}{}{}{\bfseries}{.}{ }{}

\theoremstyle{mattdef}
\newtheorem*{rmk}{Remark}

\newtheorem*{eg}{Example}

\begin{document}

\newenvironment{pfenum}{\begin{pf}\indent\begin{enumerate}\vspace{-\topsep}}{\qedhere\end{enumerate}\end{pf}}

\newcommand\ori\odot
\newcommand\tip{\triangledown}
\newcommand\ind{\operatorname{Ind}}
\newcommand\res{\operatorname{Res}}
\newcommand\iind{\operatorname{-Ind}}
\newcommand\ires{\operatorname{-Res}}
\newcommand\syd{\,\triangle\,}
\newcommand\mo{{-}1}
\newcommand\ee{\mathbf{e}}
\newcommand\reg{\mathcal{R}}
\newcommand\swed[1]{\stackrel{#1}{\wed}}
\newcommand\inv[2]{\llbracket#1,#2\rrbracket}
\newcommand\one{\mathbbm{1}}
\newcommand\bsm{\begin{smallmatrix}}
\newcommand\esm{\end{smallmatrix}}
\newcommand{\rt}[1]{\rotatebox{90}{$#1$}}
\newcommand\la\lambda
\newcommand\id{\operatorname{identity}}
\newcommand{\ol}{\overline}
\newcommand{\lan}{\langle}
\newcommand{\ran}{\rangle}
\newcommand\partn{\mathcal{P}}
\newcommand\calb{\mathcal{B}}
\newcommand\calq{\mathcal{Q}}
\newcommand\cals{\mathcal{S}}
\newcommand\calp{\mathcal{P}}
\newcommand\calf{\mathcal{F}}
\newcommand\calr{\mathcal{R}}
\newcommand\calh{\mathcal{H}}
\newcommand\calla{\mathtt{A}}
\newcommand\callb{\mathtt{B}}
\newcommand\callc{\mathtt{C}}
\newcommand{\py}[3]{\,_{#1}{#2}_{#3}}
\newcommand{\pyy}[5]{\,_{#1}{#2}_{#3}{#4}_{#5}}
\newcommand{\thmlc}[3]{\textup{\textbf{(\!\! #1 \cite[#3]{#2})}}}
\newcommand{\thmlabel}[1]{\textup{\textbf{#1}}\ }
\newcommand{\sss}{\mathfrak{S}_}
\newcommand{\weyl}{\hat{\mathfrak{S}}_}
\newcommand{\dom}{\trianglerighteqslant}
\newcommand{\doms}{\vartriangleright}
\newcommand{\ndom}{\ntrianglerighteqslant}
\newcommand{\ndoms}{\not\vartriangleright}
\newcommand{\domby}{\trianglelefteqslant}
\newcommand{\domsby}{\vartriangleleft}
\newcommand{\ndomby}{\ntrianglelefteqslant}
\newcommand{\ndomsby}{\not\vartriangleleft}
\newcommand{\subs}[1]{\subsection{#1}}
\newcommand{\subsubs}[1]{\subsubsection*{#1}}
\newcommand{\nin}{\notin}
\newcommand{\nchar}{\operatorname{char}}
\newcommand{\thmcite}[2]{\textup{\textbf{\cite[#2]{#1}}}\ }
\newcommand\zez{\mathbb{Z}/e\mathbb{Z}}
\newcommand\zsz{\mathbb{Z}/s\mathbb{Z}}
\newcommand\zepz{\mathbb{Z}/(e+1)\mathbb{Z}}
\newcommand{\bbf}{\mathbb{F}}
\newcommand{\bbg}{\mathbb{G}}
\newcommand{\bbr}{\mathbb{R}}
\newcommand{\bbc}{\mathbb{C}}
\newcommand{\bbn}{\mathbb{N}}
\newcommand{\bbq}{\mathbb{Q}}
\newcommand{\bbz}{\mathbb{Z}}
\newcommand\zo{\bbn_0}
\newcommand{\gs}{\geqslant}
\newcommand{\ls}{\leqslant}
\newcommand\dw{^\triangle}
\newcommand\wod{^\triangledown}
\newcommand{\hhh}{\mathcal{H}_}
\newcommand{\bbb}{\mathcal{B}_}
\newcommand{\aaa}{\mathcal{A}_}
\newcommand{\sect}[1]{\section{#1}}
\newcommand{\ff}{\mathfrak{f}}
\newcommand{\fff}{\mathfrak{F}}
\newcommand\cf{\mathcal{F}}
\newcommand\fkn{\mathfrak{n}}
\newcommand\fkp{\mathfrak{p}}
\newcommand\sx{x}
\newcommand\bra[1]{|#1\ran}
\newcommand\arb[1]{\widehat{\bra{#1}}}
\newcommand\foc[1]{\mathcal{F}_{#1}}
\newcommand{\clam}{\begin{description}\item[\hspace{\leftmargin}Claim.]}
\newcommand{\prof}{\item[\hspace{\leftmargin}Proof.]}
\newcommand{\malc}{\end{description}}
\newcommand\ppmod[1]{\ (\operatorname{mod}\,#1)}
\newcommand\wed\wedge
\newcommand\wede\barwedge
\newcommand\uu[1]{\,\begin{array}{|@{\,}c@{\,}|}\hline #1\\\hline\end{array}\,}
\newcommand{\ux}[1]{\operatorname{ht}_{#1}}
\newcommand\erim{\operatorname{rim}}
\newcommand\mire{\operatorname{rim}'}
\newcommand\mmod{\ \operatorname{Mod}}
\newcommand\cgs\succcurlyeq
\newcommand\cls\preccurlyeq
\newcommand\inc{\mathfrak{A}}
\newcommand\fsl{\mathfrak{sl}}
\newcommand\ba{\mathbf{s}}
\newcommand\kt[1]{|#1\rangle}
\newcommand\tk[1]{\langle#1|}
\newcommand\ket[1]{s_{#1}}
\newcommand\jn\diamond
\newcommand\UU{\mathcal{U}}
\newcommand\MM[1]{M^{\otimes#1}}
\newcommand\add{\operatorname{add}}
\newcommand\rem{\operatorname{rem}}
\newcommand\La\Lambda
\newcommand\lra\longrightarrow
\newcommand\lexg{>_{\operatorname{lex}}}
\newcommand\lexgs{\gs_{\operatorname{lex}}}
\newcommand\lexl{<_{\operatorname{lex}}}
\newcommand\tru[1]{{#1}_-}
\newcommand\ste[1]{{#1}_+}
\newcommand\out{^{\operatorname{out}}}
\newcommand\lad{\mathcal{L}}
\newcommand\hsl{\widehat{\mathfrak{sl}}}
\newcommand\fkh{\mathfrak{h}}
\newcommand\GG{H}
\newcommand\dable{restrictable}
\newcommand\infi[1]{$(\infty,#1)$-irreducible}
\newcommand\infr[1]{$(\infty,#1)$-reducible}
\newcommand\infs[1]{$#1$-signature}
\newcommand\be[2]{B^{#1}(#2)}
\newcommand\domi{dominant}
\newcommand\app{good}
\newcommand\lset[2]{\left\{\left.#1\ \middle|\ #2\right.\right\}}
\newcommand\bset[1]{\mathcal{Q}(#1)}
%abacus drawing commands - use a smallmatrix, with \bd for a bead and \nb for a space.
\newcommand\tl{\begin{picture}(8,4)
\put(4,-1){\line(0,1){3}}
\put(4,2){\line(1,0){7}}
\end{picture}}
\newcommand\tr{\begin{picture}(8,4)
\put(4,-1){\line(0,1){3}}
\put(-3,2){\line(1,0){7}}
\end{picture}}
\newcommand\tm{\begin{picture}(8,4)
\put(4,-1){\line(0,1){3}}
\put(-3,2){\line(1,0){14}}
\end{picture}}
\newcommand{\bd}{\begin{picture}(8,6)
\put(4,-1){\line(0,1){8}}
\put(4,3){\circle*{6}}
\end{picture}}
\newcommand{\nb}{\begin{picture}(8,6)
\put(4,-1){\line(0,1){8}}
\put(3,3){\line(1,0){2}}
\end{picture}}
\newcommand{\vd}{\begin{picture}(8,10)
\put(4,5){\circle*{1}}
\put(4,2){\circle*{1}}
\put(4,8){\circle*{1}}
\end{picture}}
\newcommand{\hd}{\begin{picture}(10,6)
\put(5,3){\circle*{1}}
\put(2,3){\circle*{1}}
\put(8,3){\circle*{1}}
\end{picture}}
\newcommand{\hhd}{\begin{picture}(10,12)
\put(5,6){\circle*{1}}
\put(2,6){\circle*{1}}
\put(8,6){\circle*{1}}
\end{picture}}
\newcommand{\hbd}{\begin{picture}(10,12)
\put(-2,6){\line(1,0){14}}
\put(5,6){\circle*{10}}
\end{picture}}
\newcommand{\hnb}{\begin{picture}(10,12)
\put(-2,6){\line(1,0){14}}
\put(5,4.5){\line(0,1){3}}
\end{picture}}
\newcommand{\vl}{\begin{picture}(0,6)
\put(0,-1){\line(0,1)8}\end{picture}}

\newlength\tmpl
\raggedbottom

%Topmatter
\title{The $t$-core of an $s$-core}
\author{Matthew Fayers\\\normalsize Queen Mary, University of London, Mile End Road, London E1 4NS, U.K.\\\texttt{\normalsize m.fayers@qmul.ac.uk}}
\date{}
\maketitle
\begin{center}
2000 Mathematics subject classification: 05E10, 05E18
\end{center}
\markboth{Matthew Fayers}{The $t$-core of an $s$-core}
\pagestyle{myheadings}

\begin{abstract}
We consider the $t$-core of an $s$-core partition, when $s$ and $t$ are coprime positive integers.  Olsson has shown that the $t$-core of an $s$-core is again an $s$-core, and we examine certain actions of the affine symmetric group on $s$-cores which preserve the $t$-core of an $s$-core.  Along the way, we give a new proof of Olsson's result.  We also give a new proof of a result of Vandehey, showing that there is a simultaneous $s$- and $t$-core which contains all others.
\end{abstract}

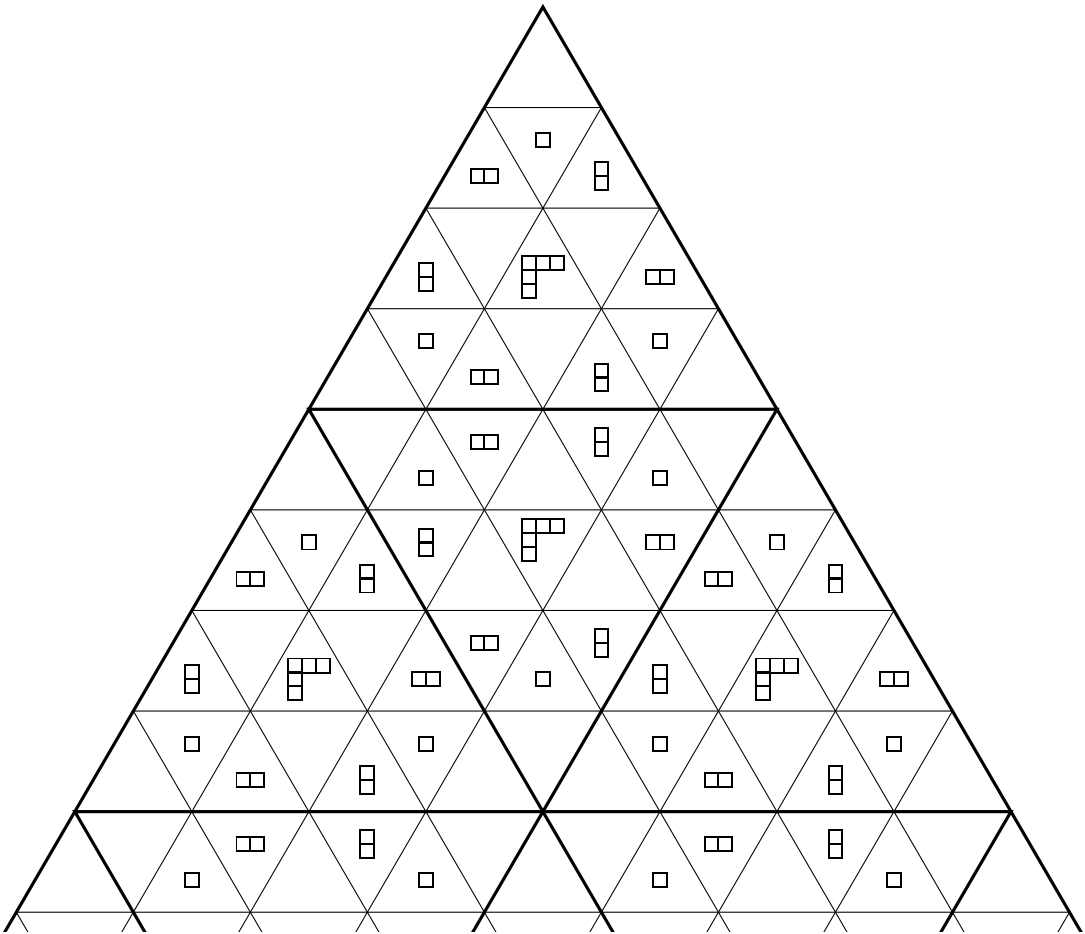
\Yboxdim{5.5pt}
\[
\begin{tikzpicture}[scale=1.71]
\clip (-1,-0.15) rectangle (9.1,8);
\draw[very thick] (-0.45,-0.78) -- (4.05,7.02) -- (8.55,-0.78);
\draw[very thick] (1.35,-0.78) -- (0.45,0.78) -- (7.65,0.78) -- (6.75,-0.78);
\draw[very thick] (3.15,-0.78) -- (5.85,3.9) -- (2.25,3.9) -- (4.95,-0.78);
\draw (0,0) -- (0.45,-0.78) -- (4.5,6.24) -- (3.6,6.24) -- (7.65,-0.78) -- (8.1,0) -- cycle;
\draw (0.9,1.56) -- (2.25,-0.78) -- (5.4,4.68) -- (2.7,4.68) -- (5.85,-0.78) -- (7.2,1.56) -- cycle;
\draw (4.95,-0.78) -- (6.75,2.34) --(1.35,2.34) -- (3.15,-0.78);
\draw (1.35,-0.78) -- (4.95,5.46) -- (3.15,5.46) -- (6.75,-0.78);
\draw (4.05,-0.78) -- (1.8,3.12) -- (6.3,3.12) -- cycle;
\draw (2.25,1.81) node{$\yng(3,1,1)$};
\draw (5.85,1.81) node{$\yng(3,1,1)$};
\draw (4.05,4.93) node{$\yng(3,1,1)$};
\draw (4.05,2.89) node{$\yng(3,1,1)$};
\foreach \x in {1.35,3.15,4.95,6.75}
\foreach \y in {0.25,1.31}
\draw (\x,\y) node{$\yng(1)$};
\foreach \x in {2.25,5.85}
\draw (\x,2.87) node {$\yng(1)$};
\foreach \x in {3.15,4.95}
\foreach \y in {3.37,4.43}
\draw (\x,\y) node{$\yng(1)$};
\foreach \y in {1.81,5.99}
\draw (4.05,\y) node{$\yng(1)$};
\foreach \x in {1.8,5.4}
\foreach \y in {0.53,1.03,2.59}
\draw (\x,\y) node{$\yng(2)$};
\foreach \x in {3.15,6.75}
\draw (\x,1.81) node{$\yng(2)$};
\foreach \y in {2.09,3.65,4.15,5.71}
\draw (3.6,\y) node{$\yng(2)$};
\foreach \x in {2.7,6.3}
\foreach \y in {0.53,1.03,2.59}
\draw (\x,\y) node{$\yng(1,1)$};
\foreach \x in {1.35,4.95}
\draw (\x,1.81) node{$\yng(1,1)$};
\foreach \y in {2.09,3.65,4.15,5.71}
\draw (4.5,\y) node{$\yng(1,1)$};
\draw (4.95,4.93) node{$\yng(2)$};
\draw (3.15,4.93) node{$\yng(1,1)$};
\draw (4.95,2.87) node{$\yng(2)$};
\draw (3.15,2.87) node{$\yng(1,1)$};

\end{tikzpicture}
\]
\Yboxdim{4.5pt}

\normalsize

\newpage

\section{Introduction}

In this paper, we define a \emph{partition} to be a infinite weakly decreasing sequence of non-negative integers with finite sum.  If $s$ is a positive integer, then a partition is an \emph{$s$-core} if it has no rim $s$-hooks, or equivalently if none of its hook lengths is divisible by $s$.  The $s$-core of an arbitrary partition is obtained by removing as many rim $s$-hooks as possible.

The notion of an $s$-core was introduced in the representation theory of the symmetric group: when $s$ is a prime, the $s$-blocks of all symmetric groups of a given defect are indexed by the $s$-cores, and the relationships between these blocks are controlled by the combinatorics of $s$-cores.  This can be generalised to representations of Iwahori--Hecke algebras of type $A$ at an $s$th root of unity (where $s\gs2$ no longer needs to be prime).

This representation-theoretic work reveals a relationship between the set of $s$-cores and the alcove geometry for the Coxeter group of type $A_{s-1}$, when $s\gs2$.  Specifically, $s$-cores are in bijection with alcoves in the dominant region of the weight space, which in turn are in bijection with cosets of the finite Coxeter group (of type $A_{s-1}$) in its affine counterpart (of type $\tilde A_{s-1}$).  Furthermore, the action of these groups of the set of alcoves can be interpreted in terms of the relationships between $s$-cores.

A recent trend in the study of $s$-cores has been to compare $s$-cores and $t$-cores, for different integers $s,t$.  For $s\gs2$ there are infinitely many $s$-cores, but if $s$ and $t$ are coprime, there are only finitely many partitions which are simultaneously $s$-cores and $t$-cores.  The exact number was found by Anderson \cite{a}, and these `$(s,t)$-cores' have since been studied in more detail.  In particular, it is known that there is an $(s,t)$-core which is larger than the others, having size $\frac1{24}(s^2-1)(t^2-1)$, and it was asked by Olsson and Stanton \cite{os} whether this partition contains all $(s,t)$-cores.  This question has been answered affirmatively by Vandehey in an unpublished thesis \cite{vand}; in the present paper we give a simpler proof.  Another avenue is pursued by Fishel and Vazirani \cite{fv}, who examine $(s,t)$-cores in connection with alcove geometry in the cases where $t\equiv \pm1\ppmod s$, exhibiting natural bijections between $(s,t)$-cores and (bounded) regions in the extended Shi arrangement.

Another aspect of the comparison of $s$- and $t$-cores is a result of Olsson \cite{ol} which says that if $s$ and $t$ are coprime and one takes the $t$-core of an $s$-core, then the resulting partition is still an $s$-core.  The main focus of this paper is to ask which $(s,t)$-core one obtains by taking the $t$-core of an $s$-core.  We explore how the symmetry of the set of $s$-cores is manifested when one replaces each $s$-core with its $t$-core.  One by-product of this is a new proof of Olsson's result.  We remark here that the hypothesis that $s$ and $t$ are coprime in Olsson's result is unnecessary, as observed by Nath \cite{na}.

Experts in combinatorial representation theory will be aware of the theory of \emph{bar partitions} and \emph{$m$-bar-cores} which control the combinatorics of spin representation of the symmetric group; it is natural to ask whether analogues of these results in the present paper hold in this context.  We address these issues in a forthcoming paper \cite{f}.

We now summarise the layout of this paper.  In Section \ref{partsec}, we give a brief account of \(s\)-cores and abacus displays.  In Section \ref{alcsec} we discuss alcove geometry and the affine Weyl group in type $A$.  We go into more detail here, since the conventions we use for alcoves are slightly unusual.  In Section \ref{mainsec} we connect $s$-cores with alcove geometry and prove our main results on the symmetry inherent in taking the $t$-core of an $s$-core.  Finally in Section \ref{lastsec} we examine the largest $(s,t)$-core, and prove that it contains all $(s,t)$-cores.

\section{Partitions}\label{partsec}

\subs{Partitions and $s$-cores}

A \emph{partition} is a sequence $\la=(\la_1,\la_2,\dots)$ of non-negative integers such that $\la_1\gs\la_2\gs\cdots$ and the sum $\la_1+\la_2+\dots$ is finite.  When writing a partition, we usually omit the trailing zeroes.  A partition is often represented by its \emph{Young diagram}, which is the set
\[[\la]=\lset{(i,j)\in\bbn^2}{j\ls \la_i}.\]
We draw the Young diagram as an array of boxes in the plane; for example, the array
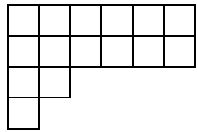
\Yboxdim{12pt}
\[\yng(6,6,2,1)\]
represents the partition $(6,6,2,1)$.  (It is usual to use a symbol such as $\varnothing$ in place of the Young diagram for the partition $(0,0,\dots)$, but in this paper we shall just use an empty diagram.)  The \emph{rim} of $\la$ is the set of nodes $(i,j)\in[\la]$ for which $(i+1,j+1)\notin[\la]$.

Now fix a positive integer $s$.  If $\la$ is a partition, a \emph{rim $s$-hook} of $\la$ is a connected portion of the rim, consisting of exactly $s$ boxes, which can be removed from $[\la]$ to leave a new Young diagram.  A partition is an \emph{$s$-core} if it has no rim $s$-hooks.  Starting from any partition $\la$ and repeatedly removing rim $s$-hooks, one obtains an $s$-core, which is independent of the choice of rim hook removed at each stage; this $s$-core is referred to as the $s$-core of $\la$.

For example, the marked boxes in the following Young diagram for $(6,6,2,1)$ constitute a rim $5$-hook.  When this is removed, the remaining partition is $(5,2,2,1)$, which has no rim $5$-hooks, and so is the $5$-core of $(6,6,2,1)$.
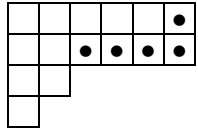
\[
\young(\ \ \ \ \ \bullet,\ \ \bullet\bullet\bullet\bullet,\ \ ,\ )
\]

The notion of an $s$-core derives from the representation theory of the symmetric group: if $s$ is a prime and $\la,\mu$ are two partitions of size $n$, the corresponding ordinary irreducible representations of $\sss n$ lie in the same $s$-block of $\sss n$ if and only if $\la$ and $\mu$ have the same $s$-core.  So the results in this paper can be interpreted as comparing the representation theory of the symmetric group for two different primes.  But from a combinatorial point of view, there is no need to assume that $s$ is prime.

\subs{Beta-numbers}\label{betanos}

Now we define beta-numbers and the abacus; these were introduced by James \cite{j0}.  Given a partition $\la$, define
\[\beta_i=\la_i-i\]
for $i\in\bbn$.  Then the sequence $\beta_1,\beta_2,\dots$ is a strictly decreasing sequence such that $\beta_i=-i$ for sufficiently large $i$.  Conversely, any such sequence is the sequence of beta-numbers of some partition.  We refer to the set $\{\beta_1,\beta_2,\dots\}$ as the \emph{beta-set} of $\la$.  Given a positive integer $s$, the \emph{$s$-runner abacus} is an abacus with $s$ infinite vertical runners, numbered $0,\dots,s-1$ from left to right; for each $j$, runner $j$ has marked positions labelled by the integers congruent to $j$ modulo $s$ increasing down the runner.  For example, the $4$-runner abacus is as follows.
\setlength\unitlength{1.1pt}
\[
\begin{smallmatrix}
0&1&2&3\\[2pt]\hline
\vd&\vd&\vd&\vd\\
\vl&\vl&\vl&\vl\\[2pt]
-4\phantom{-}&-3\phantom{-}&-2\phantom{-}&-1\phantom{-}\\[2pt]
\vl&\vl&\vl&\vl\\[2pt]
0&1&2&3\\[2pt]
\vl&\vl&\vl&\vl\\[2pt]
4&5&6&7\\[2pt]
\vl&\vl&\vl&\vl\\
\vd&\vd&\vd&\vd\\
\end{smallmatrix}
\]
The \emph{$s$-runner abacus display} for $\la$ is obtained by placing a bead on the abacus at position $\beta_i$ for each $i$.  The abacus display makes it very easy to see whether a partition is an $s$-core, since removing a rim $s$-hook corresponds to moving a bead up into an empty space immediately above it.  Hence $\la$ is an $s$-core if and only if in the $s$-runner abacus display for $\la$ every bead has a bead immediately above it.  Moreover, it is easy to obtain the abacus display for the $s$-core of $\la$ from the abacus display for $\la$: one just slides beads up until there is no bead with an empty space above it.

\begin{eg}
Suppose $\la=(6,6,2,1)$.  Then the beta-set for $\la$ is
\[\{5,4,-1,-3,-5,-6,-7,\dots\}.\]
So the $5$-runner abacus display for $\la$ is as follows.
\[
\begin{smallmatrix}
0&1&2&3&4\\[2pt]\hline
\vd&\vd&\vd&\vd&\vd\\
\bd&\bd&\bd&\bd&\bd\\
\bd&\bd&\bd&\bd&\bd\\
\bd&\nb&\bd&\nb&\bd\\
\nb&\nb&\nb&\nb&\bd\\
\bd&\nb&\nb&\nb&\nb\\
\nb&\nb&\nb&\nb&\nb\\
\nb&\nb&\nb&\nb&\nb\\
\vd&\vd&\vd&\vd&\vd
\end{smallmatrix}
\]
The abacus display of the $5$-core $(5,2,2,1)$ of $\la$ is obtained by moving the lowest bead on runner $0$ up one position.
\end{eg}

\section{Alcoves and the affine symmetric group}\label{alcsec}

In this section we introduce alcoves and the affine symmetric group.  This material will be very familiar to many readers, but we give a detailed account here because the conventions we use are slightly unusual.

\subs{Alcoves and $s$-points}

As before, we assume $s$ is a positive integer.  In fact, from now on we assume that $s\gs2$.  Let $P^s$ denote the affine space
\[P^s=\lset{p\in\bbr^s}{p_1+\dots+p_s=\mbinom s2}.\]
We define the \emph{dominant region} to be the subset of $P^s$ consisting of points $p$ for which $p_1\ls\cdots\ls p_s$.  (Note that this is unconventional -- the inequalities are usually taken the other way round.)

For each $1\ls i,j\ls s$ with $i\neq j$ and for each integer $k$, we define the hyperplane
\[
H_{ij}^k=\lset{p\in P^s}{p_j-p_i=ks},
\]
and we let $\calh=\lset{H_{ij}^k}{1\ls i<j\ls s,\ k\in\bbz}$.  The connected components of the complement in $P^s$ of the union of the hyperplanes in $\calh$ are called \emph{alcoves}.  We will abuse terminology by referring to the point $(0,1,\dots,s-1)$ as the \emph{origin}, and denoting it $\ori$.  The alcove $\calla$ containing this point is called the \emph{fundamental alcove}, and is bounded by the hyperplanes $H_{i(i+1)}^0$ (for $1\ls i<s$) and $H_{1s}^1$.

We define an \emph{$s$-point} to be a point $p\in P^s$ whose coordinates are integers which are pairwise incongruent modulo $s$.  Obviously each $s$-point is contained in some alcove, and as we shall see below, each alcove contains a unique $s$-point.

\begin{eg}
In the case $s=3$, we can draw a picture of part of the plane $P^s$ with $3$-points and hyperplanes marked as follows.
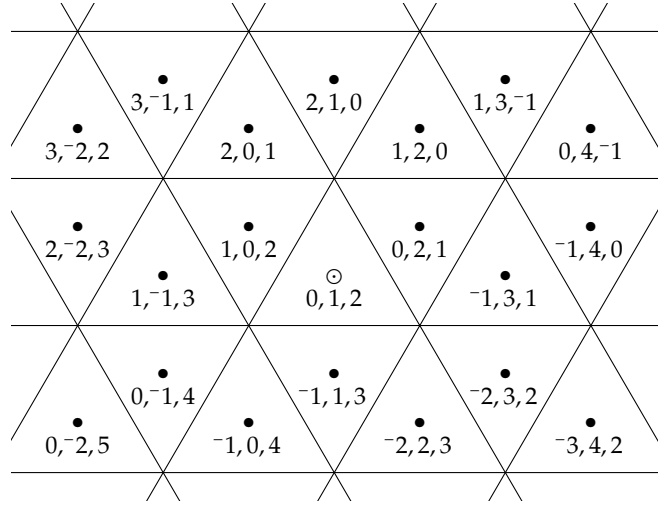
\footnotesize
\[
\begin{tikzpicture}[scale=2.5]
\clip (0.1,-0.15) rectangle (3.5,2.49);
\draw (0,0) -- (3.6,0) -- (1.8,3.12) -- cycle;
\draw (0,0.78) -- (3.6,0.78);
\draw (0,2.34) -- (3.6,2.34);
\draw (0,1.56) -- (3.6,1.56) -- (2.7,3.12) -- (0.45,-0.78);
\draw (3.15,-0.78) -- (0.9,3.12) -- (0,1.56) -- (1.35,-0.78) -- (3.6,3.12);
\draw (0,3.12) -- (2.25,-0.78) -- (3.6,1.56);
\foreach \x in {0.45,1.35,2.25,3.15}
\foreach \y in {0.26,1.82,1.3}
\draw(\x,\y)node{$\bullet$};
\foreach \x in {0.9,1.8,2.7}
\foreach \y in {0.52,2.08}
\draw(\x,\y)node{$\bullet$};
\foreach \x in {0.9,2.7}
\draw(\x,1.04)node{$\bullet$};
\draw(1.8,1.04)node{$\odot$};
\draw (0.45,0.14) node{$0,\negthinspace{}^-2,5$};
\draw (1.35,0.14) node{$\negthinspace{}^-1,0,4$};
\draw (2.25,0.14) node{$\negthinspace{}^-2,2,3$};
\draw (3.15,0.14) node{$\negthinspace{}^-3,4,2$};
\draw (0.9,0.92) node{$1,\negthinspace{}^-1,3$};
\draw (1.8,0.92) node{$0,1,2$};
\draw (2.7,0.92) node{$\negthinspace{}^-1,3,1$};
\draw (0.45,1.7) node{$3,\negthinspace{}^-2,2$};
\draw (1.35,1.7) node{$2,0,1$};
\draw (2.25,1.7) node{$1,2,0$};
\draw (3.15,1.7) node{$0,4,\negthinspace{}^-1$};
\draw (0.9,0.4) node{$0,\negthinspace{}^-1,4$};
\draw (1.8,0.4) node{$\negthinspace{}^-1,1,3$};
\draw (2.7,0.4) node{$\negthinspace{}^-2,3,2$};
\draw (0.45,1.18) node{$2,\negthinspace{}^-2,3$};
\draw (1.35,1.18) node{$1,0,2$};
\draw (2.25,1.18) node{$0,2,1$};
\draw (3.15,1.18) node{$\negthinspace{}^-1,4,0$};
\draw (0.9,1.96) node{$3,\negthinspace{}^-1,1$};
\draw (1.8,1.96) node{$2,1,0$};
\draw (2.7,1.96) node{$1,3,\negthinspace{}^-1$};
\end{tikzpicture}
\]
\normalsize
\end{eg}

For any $i,j,k$ as above, let $r_{ij}^k$ denote the orthogonal (with respect to the usual inner product on $\bbr^s$) reflection in the hyperplane $H_{ij}^k$; this is given by
\[r_{ij}^k:p\longmapsto p-(p_j-p_i-ks)(\ee_j-\ee_i),\]
where $\ee_1,\dots,\ee_s$ are the standard basis vectors.  $r_{ij}^k$ preserves the set of hyperplanes $\calh$; indeed, one can check that
\[
r_{ij}^k(H_{lm}^n)=
\begin{cases}
H_{lm}^n&(\{i,j\}\cap\{l,m\}=\emptyset)\\
H_{jm}^{n-k}&(i=l,\ j\neq m)\\
H_{lm}^{2k-n}&(i=l,\ j=m).
\end{cases}
\]
Hence the group generated by all the $r_{ij}^k$ preserves the set of alcoves.  It also preserves the set of $s$-points, and we can regard it as acting on the set of alcoves or the set of $s$-points, as appropriate.

\subs{The affine symmetric group}\label{affsymgp}

Now we consider the affine symmetric group.  This is the group $\weyl s$ with generators $\sigma_0,\dots,\sigma_{s-1}$ and relations
\begin{alignat*}2
\sigma_i^2&=1\qquad&&\text{for }i=0,\dots,s-1,\\
\sigma_i\sigma_j&=\sigma_j\sigma_i&&\text{for }i\nequiv j\pm1\ppmod s,\\
\sigma_i\sigma_j\sigma_i&=\sigma_j\sigma_i\sigma_j&\qquad&\text{for }i\equiv j+1\nequiv j-1\ppmod s.
\end{alignat*}

\subsubsection*{A level $t$ action of $\weyl s$}

There is a well-known action of $\weyl s$ on $P^s$ given by mapping the generators $\sigma_0,\dots,\sigma_{s-1}$ to the reflections in the walls of the fundamental alcove.  In fact, we give a more general version of this action.    For any positive integer $t$, define the \emph{level $t$ action} $\psi_t$ of $\weyl s$ by
\begin{alignat*}2
\sigma_i&\longmapsto r_{i(i+1)}^0\qquad&&\text{for }i=1,\dots,s-1,\\
\sigma_0&\longmapsto r_{1s}^t.
\end{alignat*}
Using the above formula for $r_{ij}^k(H_{lm}^n)$, one obtains conjugacy relations between the reflections $r_{ij}^k$, and from these it is easy to check  that this really does give an action of $\weyl s$.  Given $\sigma\in \weyl s$, we shall write $\check\sigma$ for the image of $\sigma$ under $\psi_t$, if $t$ is understood.  We may view $\psi_t$ as an action on the set of alcoves, or on the set of $s$-points, as appropriate.

It is worth while to write down explicitly the action of the generators $\sigma_0,\dots,\sigma_{s-1}$:
\begin{alignat*}2
\check\sigma_i:(p_1,\dots,p_s)&\longmapsto (p_1,\dots,p_{i-1},p_{i+1},p_i,p_{i+2},\dots,p_s)&\qquad&\text{for }i=1,\dots,s-1;\\
\check\sigma_0:(p_1,\dots,p_s)&\longmapsto (p_s-st,p_2,\dots,p_{s-1},p_1+st).
\end{alignat*}

The next lemma, which is well-known, concerns the case $t=1$.

\begin{lemma}
The image of the action $\psi_1$ includes all the reflections $r_{ij}^k$, and is transitive on the set of alcoves.
\end{lemma}

\begin{pf}
Let $G$ denote the image of $\psi_1$.  Then we have $r_{1s}^1\in G$, and we also have
\[r_{1s}^0=r_{(s-1)s}^0r_{(s-2)(s-1)}^0\dots r_{23}^0r_{12}^0r_{23}^0\dots r_{(s-2)(s-1)}^0r_{(s-1)s}^0\in G.\]
By repeatedly composing $r_{1s}^0$ and $r_{1s}^1$, we find that $r_{1s}^k\in G$ for all $k$.  And now we can see that $r_{ij}^k\in G$ for all $i,j,k$ with $i<j$ by downwards induction on $j-i$: if $j-i<s-1$, then we have either $i>1$ or $j<s$.  If $i>1$, then then by induction $r_{(i-1)j}^k\in G$, and hence $r_{ij}^k=r_{(i-1) i}^0 r_{(i-1) j}^k r_{(i-1) i}^0\in G$.  On the other hand, if $j<s$, then $r_{i(j+1)}^k\in G$ by induction, and hence $r_{ij}^k=r_{j(j+1)}^0r_{i(j+1)}^kr_{j(j+1)}^0\in G$.

To see that the action is transitive on alcoves, we note that we can get from any alcove $\callb$ to any other alcove $\callc$ by crossing some finite sequence of hyperplanes in $\calh$.  Applying the reflections in each of these hyperplanes in turn takes $\callb$ to $\callc$.
\end{pf}

Since the fundamental alcove $\calla$ clearly contains a unique $s$-point (namely the origin $\ori$), we see that each alcove contains exactly one $s$-point.  This gives a useful one-to-one correspondence between $s$-points and alcoves.

\subsubsection*{A second level $t$ action of $\weyl s$}

Now we assume that $s,t$ are coprime, and consider another level $t$ action of $\weyl s$ on the set of $s$-points.  Given an $s$-point $p$ and given $i\in\{0,\dots,s-1\}$, let $j,k$ be the unique elements of $\{1,\dots,s\}$ such that
\[p_j\equiv (i-1)t,\quad p_k\equiv it\pmod s.\]
Define $\tilde\sigma_i(p)$ to be the point obtained by replacing $p_j$ with $p_j+t$, and $p_k$ with $p_k-t$.  Clearly $\tilde\sigma_i(p)$ is an $s$-point, and it is routine to verify that the map
\[\chi_t:\sigma_i\longmapsto \tilde\sigma_i\]
defines an action of $\weyl s$ on the set of $s$-points.  Given any $\sigma\in \weyl s$, we write $\tilde\sigma$ for the image of $\sigma$ under $\chi_t$, if $t$ is understood.

Note that the maps $\tilde\sigma$ are \emph{not} isometries, and there is no natural way to extend them to give an action on the whole of $P^s$.  However, using the natural correspondence between alcoves and $s$-points, we may abuse notation and regard $\chi_t$ as an action of $\weyl s$ on the set of alcoves.  Recalling that $\calla$ denotes the alcove containing the origin $(0,1,\dots,s-1)$, we then have the following lemma, which is easy to check.

\begin{lemma}\label{agrecom}
\indent\begin{enumerate}
\vspace{-\topsep}
\item\label{commu}
If $t$ is a positive integer coprime to $s$, then the actions $\psi_t,\chi_t$ on the set of alcoves commute.
\item\label{agre}
If $t=1$ and $i\in\{0,\dots,s-1\}$, then $\check\sigma_i(\calla)=\tilde\sigma_i(\calla)$.
\end{enumerate}
\end{lemma}

Now say that two alcoves are \emph{adjacent} if there is only one hyperplane in $\calh$ separating them.

\begin{cory}
Suppose $\callb$ is an alcove, and $i\in\{0,\dots,s-1\}$, and define $\tilde\sigma_i$ using the level $1$ action $\chi_1$.  Then $\tilde\sigma_i(\callb)$ is adjacent to $\callb$.
\end{cory}

\begin{pf}
Write $\check\sigma$ for the image of $\sigma\in\weyl s$ under the level $1$ action $\psi_1$.  Since this action is transitive on the set of alcoves, we can write $\callb=\check\sigma(\calla)$ for some $\sigma$.  Hence
\begin{align*}
\tilde\sigma_i(\callb)&=\tilde\sigma_i(\check\sigma(\calla))\\
&=\check\sigma(\tilde\sigma_i(\calla))\tag*{by Lemma \ref{agrecom}(\ref{commu})}\\
&=\check\sigma(\check\sigma_i(\calla))\tag*{by Lemma \ref{agrecom}(\ref{agre})}.
\end{align*}
Clearly $\calla$ and $\check\sigma_i(\calla)$ are adjacent, and since $\check\sigma$ is an affine transformation of $P^s$, it preserves adjacency of alcoves.
\end{pf}

Using a very similar argument, one can show that if $p$ is an $s$-point and $\callb$ the alcove containing it, then each alcove adjacent to $\callb$ contains the point $\tilde\sigma_i(p)$ for some $i$.

\subsection{$s$-sets}

Define an \emph{$s$-set} to be a set of $s$ integers which are pairwise incongruent modulo $s$ and whose sum is $\mbinom s2$.  There is an $s!$-to-$1$ map from $s$-points to $s$-sets, given by forgetting the order of coordinates.  When restricted to the set of $s$-points in the dominant region, this map becomes a bijection.  Given the correspondence between $s$-points and alcoves, we have an $s!$-to-$1$ map from the set of alcoves to the set of alcoves in the dominant region; this is given by `folding' $P^s$ along the hyperplanes $H_{ij}^0$ for all $i,j$.  This folding will be useful in understanding symmetry later.

Note that our second action $\chi_t$ of $\weyl s$ on the set of $s$-points descends to an action on $s$-sets (although the action $\psi_t$ does not).  We use the same notation $\chi_t$ (and $\tilde\sigma$) for this action on $s$-sets without fear of confusion.

\section{The $t$-core of an $s$-core}\label{mainsec}

Now we come to the main part of the paper. We suppose $s,t$ are coprime integers with $s\gs2$, and we compare the $t$-cores of different $s$-cores.  By representing $s$-cores as $s$-points, we use the geometric symmetry of the last section to see the symmetry of $t$-cores of $s$-cores.

\subs{$s$-cores and $s$-sets}

Recall that a partition is an $s$-core if and only if in its $s$-runner abacus display every bead has a bead immediately above it.  This means that if we take an $s$-core $\la$ and then for each $i=0,\dots,s-1$ define $a_i$ to be the number of the highest unoccupied position on runner $i$, then $a_i\equiv i\ppmod s$ for each $i$, and
\[a_0+\dots+a_{s-1}=0+1+2+\dots+s-1=\mbinom s2.\]
Hence the set $\{a_0,\dots,a_{s-1}\}$ is an $s$-set.  We denote this $s$-set $\bset\la$, and we let $p_\la$ be the corresponding $s$-point in the dominant region, i.e.\ the point whose coordinates are the elements of $\bset\la$ arranged in ascending order.

For example, taking $s=5$ and returning to the $5$-core $(5,2,2,1)$ from the example in \S\ref{betanos}, we have
\[\bset{(5,2,2,1)}=\{5,-4,2,-2,9\},\qquad p_{(5,2,2,1)}=(-4,-2,2,5,9).\]

It is easy to check that any $s$-set is obtained from an $s$-core in this way: given an $s$-set $\calq$, construct an abacus display in which there is a bead at position $b$ if and only if there is an element of $\calq$ below $b$ on the same runner; this is then the abacus display of an $s$-core $\la$, with $\bset\la=\calq$.  Hence we have a natural bijection between $s$-cores and $s$-sets, and therefore between $s$-cores and alcoves in the dominant region.  Using this bijection, we may regard the action $\chi_t$ of the affine symmetric group on $s$-sets as an action on the set of $s$-cores.

\begin{eg}
In Figure \ref{3core} we illustrate the bijection between $3$-cores and alcoves in the dominant region of $P^3$, by drawing the Young diagram of a $3$-core inside the corresponding alcove.  (Recall that we use the empty diagram for the partition $(0,0,\dots)$.)

\Yboxdim{5pt}
\begin{figure}[ht]
\[
\begin{tikzpicture}[scale=2.4]
\clip (-1,-0.15) rectangle (4.6,3.2);
\draw (-0.45,-0.78) -- (1.8,3.12) -- (4.05,-0.78);
\draw (0,0) -- (0.45,-0.78) -- (2.25,2.34) -- (1.35,2.34) -- (3.15,-0.78) -- (3.6,0) -- cycle;
\draw (0.45,0.78) -- (1.35,-0.78) -- (2.7,1.56) -- (0.9,1.56) -- (2.25,-0.78) -- (3.15,0.78) -- cycle;
\draw (0.45,0.2) node{$\yng(6,4,2)$};
\draw (3.15,0.29) node{$\yng(3,3,2,2,1,1)$};
\draw (1.35,0.23) node{$\yng(5,3,1,1)$};
\draw (2.25,0.26) node{$\yng(4,2,2,1,1)$};
\Yboxdim{5.5pt}
\draw (0.9,0.57) node{$\yng(5,3,1)$};
\draw (1.8,0.54) node{$\yng(4,2,1,1)$};
\draw (2.7,0.51) node{$\yng(3,2,2,1,1)$};
\Yboxdim{6pt}
\draw (0.9,1.01) node{$\yng(4,2)$};
\draw (1.8,1.04) node{$\yng(3,1,1)$};
\draw (2.7,1.07) node{$\yng(2,2,1,1)$};
\Yboxdim{6.5pt}
\draw (1.35,1.33) node{$\yng(3,1)$};
\draw (2.25,1.3) node{$\yng(2,1,1)$};
\Yboxdim{7.0pt}
\draw (1.35,1.82) node{$\yng(2)$};
\draw (2.25,1.85) node{$\yng(1,1)$};
\Yboxdim{7.5pt}
\draw (1.8,2.09) node{$\yng(1)$};
\end{tikzpicture}
\]
\begin{center}
\caption{the correspondence between $3$-cores and dominant alcoves in $P^3$}
\label{3core}
\end{center}
\end{figure}
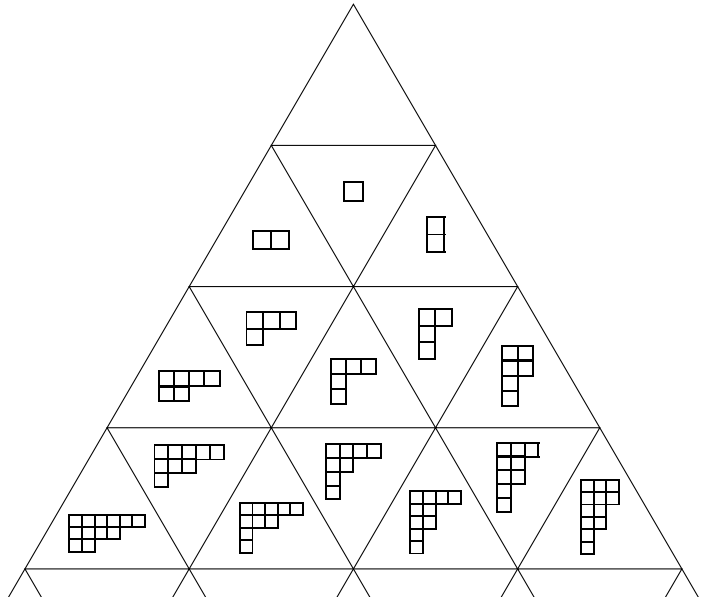

The aim of this paper is to compare the $t$-cores of different $s$-cores, when $s$ and $t$ are coprime integers.  If we take $t=4$, and expand and redraw Figure \ref{3core} with each $3$-core replaced by its $4$-core, we get the diagram on the first page of this paper.
\end{eg}

\begin{rmk}
In the case $t=1$, there is a more familiar description of the action $\chi_1$ on $s$-cores; this is described in \cite[\S4]{l} and elsewhere. Given a partition $\la$, say that a box in $[\la]$ is \emph{removable} if it can be removed to leave the Young diagram of a partition (i.e.\ it constitutes a rim $1$-hook), while a box not in $[\la]$ is an \emph{addable box of $\la$} if it can be added to $[\la]$ to give a Young diagram. The \emph{residue} of the box $(i,j)$ is defined to be the residue of $j-i$ modulo $s$.

Now when $t=1$, the action of $\tilde\sigma_k$ on an $s$-core $\la$ is given by adding all the addable boxes of residue $k$ to $\la$, or removing all the removable boxes of residue $k$ (an $s$-core cannot have addable and removable boxes of the same residue).  In terms of the abacus display for $\la$, this corresponds to interchanging the ($k-1$)th and $k$th runners of the abacus (with a slight modification in the case $k=0$).

This interpretation can be generalised to $t>1$, if one considers addable and removable rim $t$-hooks, with a suitable notion of residue.  We leave the interested reader to work out the details.
\end{rmk}

\subs{$s$-cores having the same $t$-core}

In comparing the $t$-cores of different $s$-cores, the following proposition will be crucial.

\begin{propn}\label{sametcore}
Suppose $s,t$ are coprime positive integers. Suppose $\la,\mu$ are $s$-cores, and that there is a bijection $\phi:\bset\la\to\bset\mu$ such that $\phi(b)\equiv b\ppmod t$ for every $b\in\bset\la$.  Then $\la$ and $\mu$ have the same $t$-core.
\end{propn}

\begin{pf}
We use induction on the size of $\la$. If $\la$ is not itself a $t$-core, then there is some $b$ in the beta-set for $\la$ such that $b-t$ is not in the beta-set for $\la$.  If we let $a$ be the element of $\bset\la$ which is congruent to $b$ modulo $s$, then the element of $\bset\la$ congruent to $b-t$ is $a-t-ds$ for some $d>0$.  Let $\la^-$ be the $s$-core defined by
\[\bset{\la^-}=\bset\la\cup\{a-t,a-ds\}\setminus\{a,a-t-ds\}.\]
The beta-set for $\la^-$ is obtained from the beta-set for $\la$ by replacing
\[
a-s,\ a-2s,\ \dots,\ a-ds
\]
with 
\[
a-t-s,\ a-t-2s,\ \dots,\ a-t-ds.
\]
Hence $\la^-$ is obtained from $\la$ by removing $d$ rim $t$-hooks.  So $\la^-$ is smaller than $\la$, and has the same $t$-core.  $\la^-$ and $\mu$ satisfy the hypotheses of the proposition, so by induction $\la^-$ has the same $t$-core as $\mu$, and we are done.

So we may assume that $\la$ is a $t$-core.  Symmetrically, we may assume $\mu$ is a $t$-core, and we must show that $\la=\mu$.  In other words, we need to show that a partition $\la$ which is both an $s$-core and a $t$-core is determined uniquely by the integers
\[n_i = \left|\lset{\vbox to 9pt{}a\in\bset\la}{a\equiv i\ppmod t}\right|\]
for $i=0,\dots,t-1$.

Since $s$ and $t$ are coprime, we can write the elements of $\bset\la$ as $b_0,\dots,b_{s-1}$ in such a way that $b_j\equiv -tj\ppmod s$ for each $j$.  The fact that $\la$ is a $t$-core means that for every $a$ in the beta-set for $\la$ we have $a-t$ also in the beta-set for $\la$, and this gives $b_j\gs b_{j-1}-t$ for each $j$ (reading subscripts modulo $s$).  So if we write $b_j=b_{j-1}-t+m_js$ for each $j$, then $m_0,\dots,m_{s-1}$ are non-negative integers which sum to $t$.

Now consider the following cyclic sequence of length $s+t$:
\begin{align*}
\cals=(&b_0-t,\ b_0-t+s,\ b_0-t+2s,\ \dots,\ b_1,\\
&b_1-t,\ b_1-t+s,\ b_1-t+2s,\ \dots,\ b_2,\\
&\qquad\qquad\qquad\vdots\\
&b_{s-1}-t,\ b_{s-1}-t+s,\ b_{s-1}-t+2s,\ \dots,\ b_0).
\end{align*}
The steps between consecutive terms of $\cals$ are either $-t$ ($s$ times) or $+s$ ($t$ times).  Hence modulo $t$, the steps are $0$ or $+s$.  Since $s$ generates the cyclic group of integers modulo $t$, this means that $\cals$ contains at least one term of each residue class modulo $t$.  Moreover, the terms in $\cals$ change modulo $t$ only $t$ times, so the terms in a given residue class modulo $t$ must be consecutive in $\cals$.
\clam
For each $i$, $\cals$ contains at least $n_i+1$ terms congruent to $i$ modulo $t$.
\prof
We have just seen that this is true if $n_i=0$, so suppose $n_i>0$.  $\cals$ contains all the $n_i$ elements of $\bset\la$ congruent to $i$ modulo $t$, and if we let $a$ denote the smallest of these integers, then $\cals$ also contains $a-t\notin\bset\la$.
\malc
Since $\sum_i(n_i+1)=s+t$, there must be exactly $n_i+1$ terms in $\cals$ congruent to $i$ modulo $t$ for each $i$.  And now $\cals$ is determined up to translation and cyclic permutation by the integers $n_i$: starting from the largest term divisible by $t$, $\cals$ consists of $n_0+1$ terms divisible by $t$ (with steps of $-t$ between them), and then a step of $+s$, then $n_s+1$ terms congruent to $s$ modulo $t$ (with intervening steps of $-t$), and then a step of $+s$, and so on.

So the integers $n_i$ determine $\cals$ up to translation and cyclic permutation, and hence determine $\bset\la$ up to translation.  But the sum of the elements of $\bset\la$ must be $\mbinom s2$, so $\bset\la$, and hence $\la$, is determined.
\end{pf}

Using this proposition, we see that the action $\chi_t$ of $\weyl s$ on the set of $s$-cores preserves the $t$-core of a partition.

\begin{propn}\label{orbcore}
Suppose $s$ and $t$ are coprime positive integers with $s>1$, and for $i\in\{0,\dots,s-1\}$ write $\tilde\sigma_i$ for the image of $\sigma_i$ under the action $\chi_t$.  If $\la$ is an $s$-core, then $\la$ and $\tilde\sigma_i(\la)$ have the same $t$-core.
\end{propn}

\begin{pf}
By definition, $\tilde\sigma_i$ does not change the multiset of residues modulo $t$ of the elements of an $s$-set.  Hence by Proposition \ref{sametcore}, it does not change the $t$-core of the corresponding $s$-core.
\end{pf}

We shall refer to an orbit in the set of $s$-cores under the action $\chi_t$ as a \emph{level $t$ orbit}.  From Proposition \ref{orbcore}, we see that two $s$-cores have the same $t$-core if they lie in the same level $t$ orbit.  In fact, we shall see that they have the same $t$-core only if they lie in the same level $t$ orbit; the way we do this is to show that each level $t$ orbit contains a $t$-core.  Before we do this, it will be helpful to introduce some more notation: for $s,t$ coprime, define $\calr^s_t$ to be the \emph{level $t$ rhomboid}
\[\calr^s_t=\lset{p\in P^s}{1\ls p_{i+1}-p_{i}\ls t\text{ for }i=1,\dots,s-1}.\]
For example, the rhomboid $\calr_4^3$ is the shaded portion of the dominant region of $P^3$ shown in the following diagram.
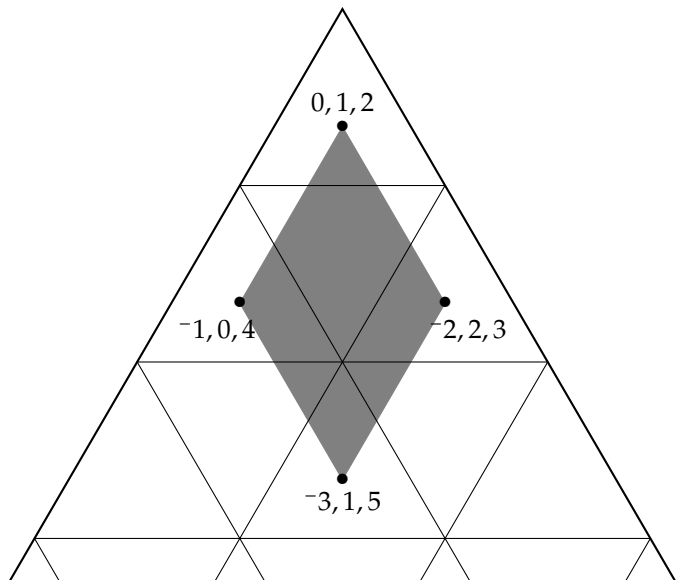
\small
\[
\begin{tikzpicture}[scale=.03]
\clip(0,60) rectangle (360,318);
\filldraw[gray] (180,104)--++(45,78)--++(-45,78)--++(-45,-78)--cycle;
\draw[thick] (0,0)--(180,312)--(360,0);
\draw (90,0)--++(-45,78)--++(270,0)--++(-45,-78)--++(-135,234)--++(90,0)--cycle;
\draw (180,0)--++(90,156)--++(-180,0)--cycle;
\draw (180,104) node{$\bullet$};
\draw (225,182) node{$\bullet$};
\draw (180,260) node{$\bullet$};
\draw (135,182) node{$\bullet$};
\draw(180,94) node{${}^-3,1,5$};
\draw(180,270) node{$0,1,2$};
\draw (125,170) node{${}^-1,0,4$};
\draw (235,170) node{${}^-2,2,3$};
\end{tikzpicture}
\]
\normalsize

Now we have the following; recall that $p_\nu$ denotes the dominant $s$-point corresponding to an $s$-core $\nu$.

\begin{propn}\label{smallest}
Suppose $s$ and $t$ are coprime integers with $s>1$, and that $O$ is a level $t$ orbit.  Let $\nu$ be an element of $O$ for which the sum $\sum_{k\in\bset{\nu}}k^2$ is minimised.  Then $p_\nu$ lies in $\calr^s_t$, and $\nu$ is a $t$-core.
\end{propn}

\begin{pf}
Suppose $p_\nu=(p_1,\dots,p_s)$ lies outside $\calr^s_t$.  By definition we have $p_1<\dots<p_s$, so the fact that $p_\nu$ is outside $\calr^s_t$ means that $p_{j+1}-p_{j}>t$ for some $j$.  If we let $f$ be the permutation of $\{1,\dots,s\}$ such that $p_{f(i)}\equiv it\ppmod s$ for each $i$, then there must be some $i$ such that $f(i-1)\ls j$ and $f(i)>j$ (reading $i-1$ modulo $s$); this then implies that $p_{f(i)}-p_{f(i-1)}>t$.   But now consider $\bset{\tilde\sigma_i(\nu)}$; this is obtained from $\bset{\nu}$ by replacing $p_{f(i-1)}$ and $p_{f(i)}$ with $p_{f(i-1)}+t$ and $p_{f(i)}-t$.  Since $p_{f(i)}-p_{f(i-1)}>t$, this gives
\[\sum_{k\in\bset{\tilde\sigma_i(\nu)}}k^2<\sum_{k\in\bset\nu}k^2,\]
a contradiction.

So $p_\nu$ lies in $\calr^s_t$.  Proving that $\nu$ is a $t$-core is very similar: if it is not, then we can find $j,k$ such that $p_k-p_j=t+as$ for some $a>0$.  Letting $i$ be such that $p_k\equiv it\ppmod s$ and applying $\tilde\sigma_i$, we derive a contradiction as above.
\end{pf}

As a consequence, we see that the element $\nu\in O$ is uniquely defined, since by Proposition \ref{orbcore} $O$ cannot contain more than one $t$-core.  Another consequence is a new proof of the following result of Olsson \cite[Theorem 1]{ol}.

\begin{thm}\label{mainolsson}
Suppose $s$ and $t$ are coprime positive integers, and $\la$ is an $s$-core.  Then the $t$-core of $\la$ is also an $s$-core.
\end{thm}

\begin{pf}
The case $s=1$ is trivial, so we may assume $s>1$.  Then by Proposition \ref{smallest}, the level $t$ orbit $O$ containing $\la$ also contains a $t$-core $\nu$.  By Proposition \ref{orbcore} $\la$ and $\nu$ have the same $t$-core, i.e.\ $\nu$ is the $t$-core of $\la$.  Since $\nu\in O$, $\nu$ is an $s$-core.
\end{pf}

Another consequence of Proposition \ref{smallest} is that two $s$-cores have the same $t$-core only if they lie in the same level $t$ orbit.

\begin{cory}
Suppose $s$ and $t$ are coprime integers, and that $\la$ and $\mu$ are $s$-cores which have the same $t$-core.  Then $\la$ and $\mu$ lie in the same level $t$ orbit.
\end{cory}

\begin{pf}
Let $\nu$ be the $t$-core of $\la$ and $\mu$.  Then $\nu$ lies in both the level $t$ orbit containing $\la$ and the level $t$ orbit containing $\mu$; so these orbits coincide.
\end{pf}

\subs{Symmetry}

Now we consider the symmetry in the diagram on the first page.  We have seen that under the action $\chi_t$ on $s$-cores, the $t$-core of an $s$-core is preserved.  However, this symmetry is obscured in the diagram on the first page because of the replacement of $s$-points by $s$-sets, or equivalently alcoves by dominant alcoves.

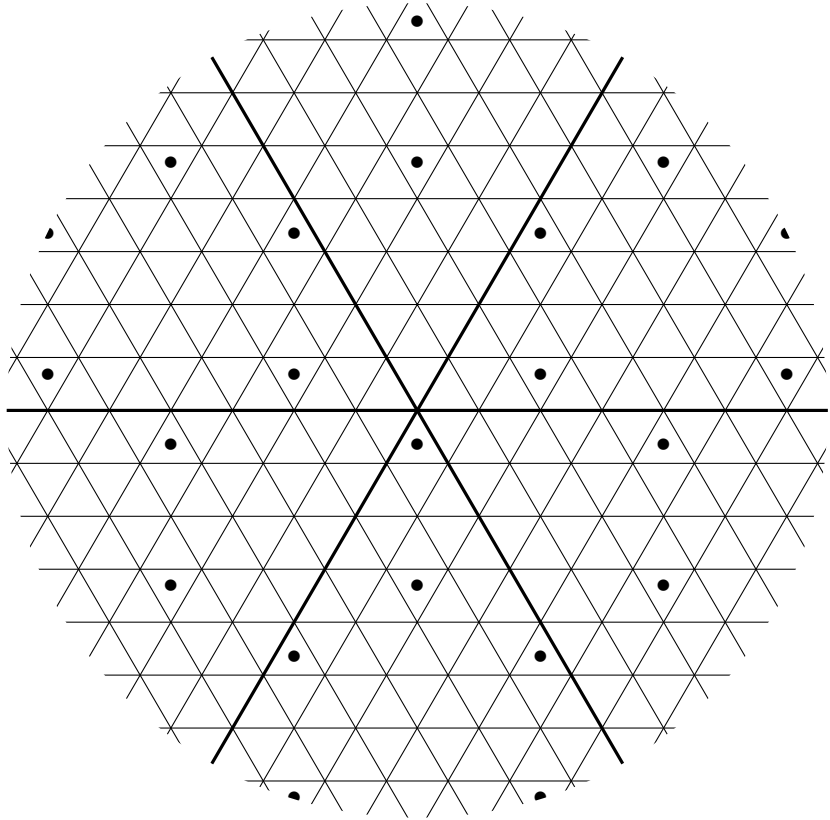
\begin{figure}[ht]
\[
\begin{tikzpicture}[scale=0.9]
\clip (0,0)circle(6);
\draw[very thick] (-3.6,-6.24) -- (3.6,6.24);
\draw[very thick] (3.6,-6.24) -- (-3.6,6.24);
\draw[very thick] (-7.2,0) -- (7.2,0);
\draw (-2.7,-6.24) -- (-6.75,0.78) -- (6.75,0.78) -- (2.7,-6.24) -- (-4.05,5.46) -- (4.05,5.46) -- cycle;
\draw (-2.7,6.24) -- (-6.75,-0.78) -- (6.75,-0.78) -- (2.7,6.24) -- (-4.05,-5.46) -- (4.05,-5.46) -- cycle;
\draw (-1.8,-6.24) -- (-6.3,1.56) -- (6.3,1.56) -- (1.8,-6.24) -- (-4.5,4.68) -- (4.5,4.68) -- cycle;
\draw (-1.8,6.24) -- (-6.3,-1.56) -- (6.3,-1.56) -- (1.8,6.24) -- (-4.5,-4.68) -- (4.5,-4.68) -- cycle;
\draw (-0.9,-6.24) -- (-5.85,2.34) -- (5.85,2.34) -- (0.9,-6.24) -- (-4.95,3.9) -- (4.95,3.9) -- cycle;
\draw (-0.9,6.24) -- (-5.85,-2.34) -- (5.85,-2.34) -- (0.9,6.24) -- (-4.95,-3.9) -- (4.95,-3.9) -- cycle;
\draw (0,-6.24) -- (-5.4,3.12) -- (5.4,3.12) -- cycle;
\draw (0,6.24) -- (-5.4,-3.12) -- (5.4,-3.12) -- cycle;
\foreach \x in {-3.6,0,3.6}
\foreach \y in {-2.6,-0.52,3.64}
\draw (\x,\y) node{$\bullet$};
\foreach \x in {-1.8,1.8}
\foreach \y in {2.6,0.52,-3.64,-5.72}
\draw (\x,\y) node{$\bullet$};
\foreach \x in {-5.4,5.4}
\foreach \y in {0.52,2.6}
\draw (\x,\y) node{$\bullet$};
\draw (0,5.72) node{$\bullet$};
\end{tikzpicture}
\]
\begin{center}
\caption{the level $4$ orbit of the origin in $P^3$}
\label{hexa}
\end{center}
\end{figure}

To show the symmetry corresponding to $\chi_t$, we consider the whole of the space $P^s$.  In our examples, we continue to take $s=3$ and $t=4$.  Figure \ref{hexa} shows part of the plane $P^3$.  The marked $3$-points are those in the level $4$ orbit containing the origin.

\begin{figure}[hbt]
\[
\begin{tikzpicture}[scale=1.1]
\clip (-4,-6) rectangle (4,0.2);
\clip (0,0)circle(6);
\draw[thick] (-3.6,-6.24) -- (0,0) -- (3.6,-6.24);
\draw (-2.7,-6.24) -- (-3.15,-5.46) -- (3.15,-5.46) -- (2.7,-6.24) -- (-0.45,-0.78) -- (0.45,-0.78) -- cycle;
\draw (-1.8,-6.24) -- (-2.7,-4.68) -- (2.7,-4.68) -- (1.8,-6.24) -- (-0.9,-1.56) -- (0.9,-1.56) -- cycle;
\draw (-0.9,-6.24) -- (-2.25,-3.9) -- (2.25,-3.9) -- (0.9,-6.24) -- (-1.35,-2.34) -- (1.35,-2.34) -- cycle;
\draw (0,-6.24) -- (-1.8,-3.12) -- (1.8,-3.12) -- cycle;
\foreach \x in {-2.25,2.25}
\draw (\x,-4.94) node{$\bullet$};
\foreach \x in {-1.35,1.35}
\foreach \y in {-4.94,-3.38,-2.86}
\draw (\x,\y) node{$\bullet$};
\foreach \x in {-1.8,1.8}
\foreach \y in {-3.64,-5.72}
\draw (\x,\y) node{$\bullet$};
\foreach \x in {-0.45,0.45}
\foreach \y in {-4.42,-1.82,-5.98}
\draw (\x,\y) node{$\bullet$};
\foreach \y in {-3.64,-2.6,-0.52,-5.72}
\draw (0,\y) node{$\bullet$};
\end{tikzpicture}
\]
\begin{center}
\caption{the level $4$ orbit of the origin in the dominant region of $P^3$}
\label{folded}
\end{center}
\end{figure}
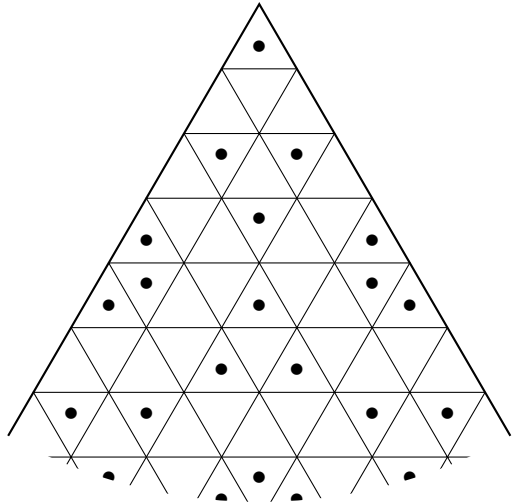

To see the corresponding orbit on $3$-sets, we fold the diagram in Figure \ref{hexa} along the bold lines (which represent the hyperplanes $H_{12}^0$, $H_{23}^0$, $H_{13}^0$).  We obtain the diagram in Figure \ref{folded}, which shows just the alcoves in the dominant region.  Comparing this with the diagram on the first page, we see that this orbit corresponds to the set of $3$-cores whose $4$-core is empty.

Now we consider the additional symmetry in the diagram on the first page: the reader will observe that the cores in the diagram are invariant under reflection in the bold lines.  These reflections (extended to the whole plane $P^3$) are the reflections contained in the image of the action $\psi_4$.  The next proposition shows that this symmetry holds in general.

\begin{propn}
Suppose $s,t$ are coprime, and $p,q$ are $s$-points which lie in the same orbit under the level $t$ action $\psi_t$.  Then the $s$-cores corresponding to $p,q$ have the same $t$-core.
\end{propn}

\begin{pf}
This is immediate from Proposition \ref{sametcore} and the formul\ae{} for $\check\sigma_0,\dots,\check\sigma_{s-1}$ in \S\ref{affsymgp}.  
\end{pf}

Note that, unlike the orbits for the action $\chi_t$, different orbits under $\psi_t$ can yield the same $t$-core; this can be seen in the diagram on the first page, where there are six different orbits yielding the empty partition.

There is further symmetry in the diagram on the first page.  The image of the action $\psi_t$ is generated by the reflections in the hyperplanes $H_{i(i+1)}^0$ for $1\ls i<s$, and $H_{1s}^t$.  These hyperplanes bound a simplex which may be obtained by dilating the fundamental alcove by a factor of $t$, fixing the point $x_0=(\frac{s-1}2,\dots,\frac{s-1}2)$.  The vertices of this simplex are the points $x_0,\dots,x_{s-1}$ defined by
\[
(x_i)_j=
\begin{cases}
\mfrac{s-1}2+(i-s)t&(j\ls i)\\[6pt]
\mfrac{s-1}2+it&(j>i).
\end{cases}
\]
There is a unique affine-linear map $\alpha_t:P^s\to P^s$ which permutes these vertices cyclically: this is given by
\[\alpha_t(p_1,\dots,p_s)=(p_s-(s-1)t,p_1+t,p_2+t,\dots,p_{s-1}+t).\]
One can check that this map preserves the set of $s$-points, and also the set of alcoves.  Using Proposition \ref{sametcore} and the formula above, we see that if $p$ is an $s$-point, then the $s$-cores corresponding to $p$ and $\alpha_t(p)$ have the same $t$-core.  In the diagram on the first page, this can be seen as a rotational symmetry: the dilated fundamental alcove is the large triangle at the top of the diagram bounded by bold lines, and one can see that rotating this triangle through angle $\frac{2\pi}3$ preserves the $t$-cores in the diagram.

\section{The largest $(s,t)$-core}\label{lastsec}

We continue to assume that $s,t$ are coprime positive integers with $s\gs2$.  An \emph{$(s,t)$-core} means a partition which is both an $s$-core and a $t$-core.  It can be inferred from the results in previous sections that there are only finitely many $(s,t)$-cores.  In fact, the exact number of $(s,t)$-cores is $\mfrac1{s+t}\mbinom{s+t}s$, as was shown by Anderson \cite[Theorems 1 \& 3]{a}.

The aim in this section is to show that one of these $(s,t)$-cores contains all the others.  To describe this $(s,t)$-core and to set up the proof, recall the level $t$ rhomboid
\[\calr^s_t=\lset{p\in P^s}{1\ls p_{i+1}-p_{i}\ls t\text{ for }i=1,\dots,s-1}.\]
Proposition \ref{smallest} implies that if $\nu$ is an $(s,t)$-core, then the point $p_\nu$ lies in $\calr^s_t$.

The vertex of $\calr^s_t$ opposite the origin will be denoted $\tip$; this has coordinates
\[\left(\frac{s-1+t(1-s)}2,\frac{s-1+t(3-s)}2,\frac{s-1+t(5-s)}2,\dots,\frac{s-1+t(s-1)}2\right).\]
Since at least one of $s,t$ is odd, these coordinates are integers, and so $\tip$ is an $s$-point.  Looking at the beta-set of the corresponding $s$-core, we see that it is a $t$-core (cf.\ the proof of Proposition \ref{sametcore}).  Following \cite{os} we denote this partition $\kappa_{s,t}$.  This partition has been studied before; it was shown by Kane \cite{k} that its size is $\frac1{24}(s^2-1)(t^2-1)$, and this is known to be the largest size of any $(s,t)$-core.  In fact, the following stronger statement is true.

\begin{thm}\thmcite{vand}{Theorem 2.4}\label{maincontain}
Suppose $s$ and $t$ are coprime positive integers, and $\la$ is an $(s,t)$-core.  Then $[\la]\subseteq[\kappa_{s,t}]$.
\end{thm}

This theorem answers a question of Olsson and Stanton \cite[Remark 4.11]{os}, who prove the theorem in the case $t=s+1$.  Our aim in this section is to give a new proof of this theorem.

We shall prove Theorem \ref{maincontain} using the correspondence between alcoves and $s$-cores.  First we need a result about hyperplanes meeting $\calr^s_t$.

\begin{lemma}\label{meetr}
Suppose $1\ls i<j\ls s$ and $k\in\bbz$ are such that the hyperplane $H_{ij}^k$ intersects $\calr^s_t$.  Then the origin and the point $\tip$ lie on opposite sides of $H_{ij}^k$.
\end{lemma}

\begin{pf}
For every point in $p\in\calr^s_t$ we have $j-i\ls p_j-p_i\ls (j-i)t$, so the statement that $H_{ij}^k$ meets $\calr^s_t$ implies that $\mfrac{j-i}s<k<\mfrac{(j-i)t}s$ (equality cannot hold on either side because $\mfrac{j-i}s,\mfrac{(j-i)t}s$ are not integers).  Writing $\ori=(\ori_1,\dots,\ori_s)$ for the origin, we have $\ori_j-\ori_i=j-i$; on the other hand, writing $\tip=(\tip_1,\dots,\tip_s)$, we have $\tip_j-\tip_i=(j-i)t$.  Hence $\ori$ and $\tip$ lie on opposite sides of $H_{ij}^k$.
\end{pf}

For the rest of this section $\tilde\sigma_0,\dots,\tilde\sigma_{s-1}$ denote the images of the generators $\sigma_0,\dots,\sigma_{s-1}$ \emph{under the level $1$ action} $\chi_1$ (not the level $t$ action) of $\weyl s$.

First we show that if $\la$ is an $s$-core, then one of the Young diagrams $[\la]$, $[\tilde\sigma_i(\la)]$ contains the other.

\begin{lemma}\label{contains}
Suppose $\la$ is an $s$-core, and $i\in\{0,\dots,s-1\}$.  Let $a$ and $b$ be the elements of $\bset\la$ congruent to $i-1$ and $i$ modulo $s$, respectively.  
If $b\ls a+1$, then $[\tilde\sigma_i(\la)]\supseteq[\la]$.
\end{lemma}

\begin{pf}
We obtain $\bset{\tilde\sigma_i(\la)}$ from $\bset\la$ by replacing $a,b$ with $a+1,b-1$.  Since $b\ls a+1$, this means that the beta-set for $\tilde\sigma_i(\la)$ is obtained from the beta-set for $\la$ by increasing each of the numbers
\[
b-1,\,b-1+s,\,b-1+2s,\,\dots,\,a-s
\]
by $1$.  So for each $j$ the $j$th beta-number for $\tilde\sigma_i(\la)$ is either equal to or one greater than the $j$th beta-number for $\la$; hence for each $j$ we have either $(\tilde\sigma_i(\la))_j=\la_j$ or  $(\tilde\sigma_i(\la))_j=\la_j+1$.  And so $[\tilde\sigma_i(\la)]\supseteq[\la]$.
\end{pf}

We now consider this containment relation in terms of alcoves.  Recall that for any alcove $\callb$ and for any $i$, the alcove $\tilde\sigma_i(\callb)$ is adjacent to $\callb$.

\begin{lemma}\label{sepalc}
Suppose $p$ is an $s$-point, and $i\in\{0,\dots,s-1\}$.  Let $\callb$ be the alcove containing $p$, let $H$ be the unique hyperplane in $\calh$ separating $\callb$ and $\tilde\sigma_i(\callb)$, and let $\la$ be the $s$-core corresponding to $p$.  Suppose the origin lies on the same side of $H$ as $p$.  Then $[\tilde\sigma_i(\la)]\supseteq[\la]$.
\end{lemma}

\begin{pf}
Write $p=(p_1,\dots,p_s)$, and let $j,l$ be such that $p_j\equiv i-1$ and $p_l\equiv i\ppmod s$.  Then $\tilde\sigma_i(p)$ is obtained by replacing $p_j$ with $p_j+1$ and $p_l$ with $p_l-1$.

If we let $k=\mfrac{p_j-p_l+1}s$, then we see that $p$ and $\tilde\sigma_i(p)$ lie on opposite sides of $H_{lj}^k$.  Hence we must have $H=H_{lj}^k$.  We have $p_j-p_l<sk$, so the fact that the origin lies on the same side of $H$ as $p$ means that $j-l<sk$, whence $k\gs0$, so that $p_l\ls p_j+1$.  Now the result follows from Lemma \ref{contains}.
\end{pf}

Now we can prove Theorem \ref{maincontain}.  In fact, we prove a stronger statement.

\begin{propn}\label{stronger}
Suppose $s,t$ are coprime and $p$ is an $s$-point in $\calr^s_t$, and let $\la$ be the corresponding $s$-core.  Then $[\la]\subseteq[\kappa_{s,t}]$.
\end{propn}

\begin{pf}
Consider the hyperplanes in $\calh$ which meet the line segment joining $p$ to $\tip$.  Each of these hyperplanes meets $\calr^s_t$, so by Lemma \ref{meetr} has the origin and $\tip$ on opposite sides of it.

By slightly deforming the line segment from $p$ to $\tip$ if necessary, we can construct a path which meets each of these hyperplanes once without meeting two of them simultaneously, and does not meet any other hyperplane in $\calh$.  If we let $\callb_0,\callb_1,\dots,\callb_d$ be the alcoves meeting this path, and $p^0,\dots,p^d$ the $s$-points contained in these alcoves, then we have $p^0=p$ and $p^d=\tip$.  Moreover, for each $l=1,\dots,d$, the points $p^{l-1}$ and $p^l$ lie in adjacent alcoves, so there is some $i_l$ such that $p^l=\tilde\sigma_{i_l}(p^{l-1})$.  If we let $H_{(l)}$ be the hyperplane separating $\callb_{l-1}$ and $\callb_l$, then $p^l$ lies on the same side of $H_{(l)}$ as $\tip$, and so the origin lies on the opposite side of $H_{(l)}$ to $p^l$. Hence, if we write $\la^{(l)}$ for the $s$-core corresponding to $p^l$ for each $l$, then by Lemma \ref{sepalc} we have
\[
[\la]=[\la^{(0)}]\subseteq[\la^{(1)}]\subseteq\dots\subseteq[\la^{(d)}]=[\kappa_{s,t}].\qedhere
\]
\end{pf}

Since $p_\nu\in\calr^s_t$ for every $(s,t)$-core $\nu$, Proposition \ref{stronger} implies Theorem \ref{maincontain}.

\begin{rmk}
In fact, our results can be interpreted as a rather stronger result than just saying that $\kappa_{s,t}$ contains every $(s,t)$-core, by considering further the Coxeter group $\weyl s$.  We summarise this very briefly.  The level $1$ action of $\weyl s$ on alcoves is faithful; under the $s!$-to-$1$ map that sends alcoves to dominant alcoves, the alcoves that map to the fundamental alcove $\calla$ are precisely those of the form $\tilde\sigma(\calla)$ for $\sigma$ in the parabolic subgroup $\sss s$ of $\weyl s$ generated by $\sigma_1,\dots,\sigma_{s-1}$.  Hence there is a bijection between dominant alcoves and left cosets of $\sss s$ in $\weyl s$, given by $\tilde\sigma(\calla)\leftrightarrow\sigma\sss s$.  Given our bijection between dominant alcoves and $s$-cores, we have a bijection between the set of $s$-cores and the set $\weyl s/\sss s$ of left cosets of $\sss s$ in $\weyl s$.

There are two well-known partial orders on $\weyl s/\sss s$: the \emph{Bruhat order} and the \emph{left order}.  To describe these, we need to recall that each left coset contains a unique \emph{minimal representative}, i.e.\ an element whose length with respect to the generators $\sigma_0,\dots,\sigma_{s-1}$ is minimised.  Given two elements $\sigma,\tau$ which are minimal representatives of their cosets $\sigma\sss s$ and $\tau\sss s$, we say that $\sigma\sss s\ls\tau\sss s$ in the Bruhat order if there is a reduced expression for $\tau$ from which we can delete some terms to yield an expression for $\sigma$.  On the other hand, we say that $\sigma\sss s\ls \tau\sss s$ in the left order if there is a reduced expression for $\tau$ from which we can delete an initial segment to yield an expression for $\sigma$.  Obviously the Bruhat order is a refinement of the left order.

Now suppose $\la$ and $\mu$ are the $s$-cores corresponding to $\sigma,\tau$.  It is shown in \cite[Proposition 4.1]{l} that $[\la]\subseteq[\mu]$ if and only if $\sigma\sss s\ls \tau\sss s$ in the Bruhat order.  What the results in this section show is that if $\nu$ is an $(s,t)$-core and $\pi\sss s$, $\rho\sss s$ are the cosets corresponding to $\nu$ and $\kappa_{s,t}$, then $\pi\sss s\ls\rho\sss s$ in the left order, which is a stronger condition.
\end{rmk}

\end{document}